\documentstyle[12pt]{article}
\newcommand{\C}{{\bf C}}
\newcommand{\Z}{{\bf Z}}
\newtheorem{thm}{Theorem}

\begin{document}
\title{Quantum Galois theory for finite groups}
\author{Akihide Hanaki and Masahiko Miyamoto \\
        and Daisuke Tambara}
\date{}
\maketitle
 
Dong and Mason \cite{DM1} initiated a systematic research
for a vertex operator algebra with a finite automorphism group,
which is referred to as the ``operator content of orbifold models''
by physicists \cite{DVVV}.
The purpose of this paper is to extend one of their main results.
We will assume that the reader is familiar with
the vertex operator algebras (VOA), see \cite{B},\cite{FLM}.
 
Throughout this paper, $V$ denotes a simple vertex operator algebra,
$G$ is a finite automorphism group of $V$,
$\C$ denotes the complex number field,
and $\Z$ denotes rational integers.
Let $H$ be a subgroup
of $G$ and $Irr(G)$ denote the set of all irreducible
${\C}G$-characters.
In their paper \cite{DM1}, they studied the sub VOA
$V^H=\{v\in V: h(v)=v \mbox{\ for all\ } h\in H\}$ of $H$-invariants and
the subspace $V^{\chi}$ on which
$G$ acts according to $\chi\in Irr(G)$. Especially,
they conjectured the following
Galois correspondence between sub VOAs of $V$ and subgroups of $G$ and
proved it for an Abelian or dihedral group $G$ \cite[Theorem 1]{DM1} and 
later for nilpotent groups \cite{DM2},
which is an origin of their title of \cite{DM1}. 
 
\bigskip
\noindent
{\bf Conjecture} {(\bf Quantum Galois Theory)} \
Let $V$ be a simple VOA and $G$ a finite and faithful group of
automorphisms of $V$. Then there is a bijection between the subgroups
of $G$ and the sub VOAs of $V$ which contains $V^G$ defined by the map
$H \to V^H$.
\bigskip
 
Our purpose in this paper is to prove the above conjecture.   
Namely, we will prove:
 
\begin{thm}
Let $V$ be a simple VOA and $G$ a finite and faithful group of
automorphisms of $V$.
Then there is a bijection between the subgroups
of $G$ and the sub VOAs of $V$ which contains $V^G$ defined by the map
$H \ (\leq G) \to V^H \ (\supseteq V^G)$.
\end{thm}
 
We adopt the notation and results in \cite{DM1} and \cite{DLM}.
Especially, the following result in \cite{DLM}
is the main tool for our study.
 
\begin{thm}[DLM,Corollary 2.5]
Suppose that $V$ is a simple VOA and that $G$ is a finite and faithful
group of automorphisms of $V$.  Then the following hold: \\
(i) For $\chi\in Irr(G)$, each $V^{\chi}$ is a simple module for
the $G$-graded VOA  ${\C}G\otimes V^G$ of the form
$$  V^{\chi}=M_{\chi}\otimes V_{\chi}   $$
where $M_{\chi}$ is the simple ${\C}G$-module affording $\chi$ and where
$V_{\chi}$ is a simple $V^G$-module.  \\
(ii)  The map $M_{\chi} \to V_{\chi}$ is a bijection from the set of
simple ${\C}G$-modules to the set of inequivalent simple $V^G$-modules
which are contained in $V$.
\end{thm}
 
It was proved in \cite[Lemma 3.2]{DM1} that
the map $H\ (\leq G) \to V^H\ (\supseteq V^G)$ is
injective.  Therefore, it is sufficient to show that
for any sub VOA $W$ containing $V^G$, there is a subgroup $H$ of $G$
such that $W=V^H$.
Our first purpose is to transform the assumption
of quantum Galois
conjecture to the following purely group theoretic condition :
 
\bigskip
\noindent
{\bf Hypotheses (A)} \
Let $G$ be a finite group and $\{M_{\chi}:\chi\in Irr(G)\}$
be the set of all nonisomorphic simple modules of $G$.
Assume $M_{1_G}$ is a trivial module.
Let $R$ be a subspace of $M=\oplus_{\chi\in Irr(G)} M_{\chi}$
containing $M_{1_G}$.
Assume that for any $G$-homomorphism $\pi:M\otimes M \to M$,
$\pi(R\otimes R)\subseteq R$.
\bigskip
 
Let's show how to transform the assumption of quantum Galois conjecture
to Hypotheses (A).
We first explain the way to give a relation between $V$ and $M$.
Let introduce a relation
$u\sim v$ for two nonzero elements $u,v\in V$
if there is an element $w\in V^G$ such that
$v=w_nu$ or $u=w_nv$ for some $n\in Z$.  Extend this relation
into an equivalent relation $\equiv $ as follows: \\
$u\equiv v$ if there are $u^1,...,u^m\in V$ such that
$u\sim u^1\sim...\sim u^m\sim v$. \\
This equivalent relation implies that
for any $G$-homomorphism $\phi:V \to M_{\chi}$, the image of each equivalent
class is uniquely determined up to scalar times. Also since $V_{\chi}$
is spanned by $\{ v_ns_{\chi}: v\in V^G, n\in \Z \}$ for some
$s_{\chi}\in V_{\chi}$ by
$[$DM1 Proposition 4.1$]$, $\phi(V_{\chi})={\C}\phi(s_{\chi})$ is a subspace
of dimension at most one.
 
Let's start the transformation. 
We recall $$ V=\oplus_{\chi} (M_{\chi}\otimes V_{\chi}) $$
from Theorem 2 ([DLM, Corollary 2.5]). We may view $M_{\chi}$ as 
a subspace of $V$. Set $M=\oplus M_{\chi}$. 
 Let $W$ be a sub VOA containing $V^G$.
Let $U_{\chi}$ be the $V^G$-subspace of $W$ whose composition factors are
all isomorphic to $V_{\chi}$.  Then
$U_{\chi}=W\cap (M_{\chi}\otimes V_{\chi})$ 
and so $W=\oplus_{\chi} U_{\chi}$.  Set
$$R_{\chi}=\{m\in M_{\chi}| m\otimes V_{\chi}\subseteq W \}$$
and $R=\oplus R_{\chi}$. In particular, we have 
$$ W=\oplus_{\chi} (R_{\chi}\otimes V_{\chi}).$$
Replacing $M_{\chi}$ by $L(k_1)...L(k_r)M_{\chi}\cong M_{\chi}$ 
if necessary, we may think that 
$R=\oplus R_{\chi}$ and $M=\oplus M_{\chi}$ are 
subspaces of a homogeneous part $V_p$ of $V$ for some $p$. 
Let $\{v^1,...,v^s\}$ be a basis of $R$ and
$\{v^1,...,v^s,...,v^n\}$ 
be a basis of $M$.  The proof of Lemma 3.1 in \cite{DM1} shows that 
$\{ Y(v^i,z)v^j\ :\ i,j=1,...,n \}$ is a linearly independent set.  
Define $\pi_m: M\times M \to V_m$ by $\pi_m(v^i\times v^j)=(v^i)_mv^j$.  
Since $\{ Y(v^i,z)v^j \}$ is a linearly independent set, 
$\cap_{m\in {\bf Z}}{\rm Ker}(\pi_m)=0$.  
Since $M\times M$ has only a finite 
dimension, there is a finite set $\{ a,a+1,...,b\}$ of integers such 
that $\cap_{m=a}^b{\rm Ker}(\pi_m)=0$.  
Since  the grade of $(v^i)_m(v^j)$ depends only on $m$ and 
$(v^i)_m(v^j)$ and $(v^h)_{m'}(v^k)$ belong to different 
homogeneous spaces for $m\not= m'$, 
the map $\pi: M\times M \to V$ given by 
$\pi(v^i\times v^j)=\sum_{m=a}^b\{(v^i)_m(v^j)\}$
is injective.  Set 
$E={\rm Im}(\pi)=<\sum_{m=a}^b (v^i)_m(v^j): i,j=1,...,n>$. 
Clearly, $E$ is a $G$-invariant subspace of $V$. Decompose $V$ into 
a direct sum $V=E\oplus E'$ of $E$ and some $G$-submodule $E'$ of $V$.  
Define $\mu:V=E\oplus E'\to M\times M$ by 
$\mu(e+e')=\pi^{-1}(e)$ for $e\in E, e'\in E'$.   
This is a $G$-epimorphism. 
Since $W$ is a sub VOA, the any products $u_nv$ of two elements $u,v$
of $W$ are in $W$.  Hence, $\pi(R\times R) \subseteq W$ and so 
the image $\mu(W)$ contains $R\times R$.  
Therefore, 
for any $G$-homomorphism $\phi:M\otimes M\to M$, we have  
$\phi(R\times R)\subseteq \phi(\mu(W))$. 
On the other hand, for any $G$-homomorphism $\psi:V \to M$, 
we have $\psi(W)\subseteq R$ since 
$W=\oplus_{\chi} (R_{\chi}\otimes V_{\chi})$.
Hence, we have 
$$\phi(R\times R)\subseteq \phi(\mu(W))=(\phi\mu)(W)\subseteq R$$ 
for any $G$-homomorphism $\phi:M\otimes M\to M$. 
Namely, $R$ satisfies the Hypotheses (A). \\
 
We will next prove the following group theoretic problem:
 
\begin{thm}
 Let $G$ be a finite group and $\{M_{\chi}:\chi\in Irr(G)\}$
be the set of all simple modules of $G$.  Assume $M_1={\C}$
is a trivial module.
Let $R$ be a subspace of $M=\oplus_{\chi\in Irr(G)} M_{\chi}$
containing $M_1$.
Assume that $R$ satisfies the following condition:
for any $G$-homomorphism $\pi:M\otimes M \to M$,
$\pi(R\otimes R)\subseteq R$.
Then there is a subgroup $G_1$ of $G$ such that $R=M^{G_1}$.
\end{thm}
 
{\it Proof.}
Consider the group algebra $\C G$,
and write
$\C G=\oplus_{\chi\in Irr(G)}M_\chi{}^{\chi(1)}$ as a left
$G$-module.
We define a subspace $S$ of $\C G$ by
$$S=\bigoplus_{\chi\in Irr(G)} (R\cap M_\chi)^{\chi(1)}.$$
Then $S$ is a right ideal of $\C G$.
Note that $R=\oplus_\chi (R\cap M_\chi)$, since $M_1={\C}\subseteq R$
and $\pi_{\chi}({\C}\otimes R)\subseteq R$ for a projection
$\pi_{\chi}:{\C}\otimes M \to M_{\chi}\subseteq M$.
Thus  $S$ satisfies that
$\pi'(S\otimes S)\subseteq S$ for any $G$-homomorphism $\pi':\C
G\otimes\C G \to\C G$.
 
We define a new product $\circ$ in $\C G$ by
$$\left(\sum_{g\in G} a_g g\right)\circ\left(\sum_{g\in G} b_g g\right)
  = \sum_{g\in G} a_g b_g g,$$
where $a_g$, $b_g\in \C$.
Then $(\C G, \circ)$ is a semisimple commutative algebra
with the identity $\sum_{g\in G}g$.
We write the identity by $1^\circ$.
 
Define $\pi':\C G\otimes\C G\to \C G$ by
$$\pi'\left(\left(\sum_{g\in G}a_g g\right)
  \otimes\left(\sum_{g\in G}b_g g\right)\right)
  = \sum_{g\in G} a_g b_g g.$$
Then $\pi'$ is a $G$-homomorphism,
and so $\pi'(S\otimes S)\subseteq S$.
This means that $S$ is a subalgebra of $(\C G,\circ)$.
 
Since $R$ contains the trivial module,
$1^\circ$ belongs to $S$.
Consider the primitive idempotent decomposition of $1^\circ$ in
$S$.
Then there exists a partition $G=\cup_i G_i$
such that $\sum_{g\in G_i}g$ is a primitive idempotent in $S$
for any $i$.
Put $e_i=\sum_{g\in G_i}g$,
then $S=\oplus_i \C e_i$ since $S$ is semisimple and commutative.
Assume $1_G\in G_1$.
For $h\in G_1$,
$e_1 h$ is in $S$ since $S$ is a right ideal in $\C G$.
By the form of $S$, $e_1 h$ is a sum of some $e_i$'s.
But $h\in G_1h$, so $e_1h=e_1$.
This means $G_1$ is a subgroup of $G$.
Similarly $G_i$ is a left coset of $G_1$ in $G$ for any $i$.
 
Now $S=\C G^{G_1}$, and so $R=V^{G_1}$.
The proof is completed.
\hfill$\Box$
\bigskip
 
Let's go back to the proof of Theorem 1 (the quantum Galois theory).
By the above theorem 3, there is a subgroup $H$ of $G$ such that
$R_{\chi}=M_{\chi}^H$ and so $$U_{\chi}=M_{\chi}^H\otimes V_{\chi}
=(M_{\chi}\otimes V_{\chi})^H.$$  Hence, $W=V^H$.
This completes the proof of Theorem 1.
 
\bigskip
{\bf Acknowledgment} \\
We would like to express our thanks to C.Lam, 
C.Dong and G.Mason for their helpful comments.

\medskip
{\small\sc
\noindent
Akihide Hanaki \\
Faculty of Engineering,
Yamanashi University,
Kofu 400, Japan \\
 
\noindent
Masahiko Miyamoto \\
Institute of Mathematics,
University of Tsukuba,
Tsukuba 305, Japan \\
 
\noindent
Daisuke Tambara\\
Department of Mathematics,
Hirosaki University,
Hirosaki 036, Japan
}
\end{document}